\documentclass[12pt,epic,eepic] {article}
\usepackage{amssymb}
\usepackage{setspace}
\usepackage{graphicx}
\newcommand{\R}{\mathbb R}
\newcommand{\N}{\mathbb N}

\begin{document}
\onehalfspacing
\title{On a question by Corson about point-finite coverings}
\author{A. Marchese and C. Zanco}
\date{}
\maketitle
\bigskip
To appear in Israel J. Math.

\footnotetext{Research of the first author was supported
by the GNAMPA of the Istituto Nazionale di Alta Matematica of Italy;
research of
the second author was supported in part by the GNAMPA of the
Istituto Nazionale di Alta Matematica of Italy and in part by the Center
for Advanced Studies in Mathematics at the
Ben-Gurion University of the Negev, Beer-Sheva, Israel.}

\bigskip
\bigskip
\bigskip
\bigskip

\begin{abstract}
We answer in the affirmative the following question raised by H. H. Corson in 1961: " Is it
possible to cover every Banach space X by bounded convex sets with nonempty interior in such a way
that no point of X belongs to infinitely many of them?"

Actually we show the way to produce in every Banach space X a bounded convex tiling of order 2,
i.e. a covering of X by bounded convex closed sets with nonempty interior (tiles)
such that the interiors are pairwise disjoint and no point of X belongs to more than two tiles.

\end{abstract}

\bigskip

\bigskip

\bigskip

{\it 2000 Mathematics Subject Classification}: 46B20, 52A45.

\bigskip

 {\it Key words and phrases}: coverings of Banach spaces, point-finite coverings, tilings.

\vfill\eject

\section{Introduction, notation, main statement}

\ \ \ Throughout the paper, by {\it covering} of a Banach space
$X$ we mean a family $\{A_\lambda\}_{\lambda \in \Lambda}$ of
proper subsets of $X$ ($\Lambda$ any set of indices) such that $X
= \cup_\lambda A_\lambda$. By {\it body} in $X$, we mean a
nonempty proper subset of $X$ that is
contained in the closure of its connected interior.
A covering of $X$ is called a {\it tiling} of $X$ whenever its
members ({\it tiles}) are closed bodies with pairwise disjoint
interiors. A covering $\tau$ of $X$ is said to be {\it
point-finite} if no point of $X$ belongs to infinitely many members of
$\tau$; the (possibly infinite) {\it order} of $\tau$ is the
supremum of those $n$ in ${\N}$ such that there exist $n$ members
of $\tau$ with a common point. A covering $\tau$ of $X$ is said to
be {\it locally finite} if every point in $X$ has a neighborhood
that meets only finitely many members of $\tau$ (equivalently, if
every compact subset of $X$ meets only finitely many members of
$\tau$). It is easy to produce examples of point-finite coverings
(even by convex bodies, even tilings) that are not locally finite: for an
exhaustive discussion of this topic in the general setting of spaces of any
dimension see the nice paper \cite{K1} by V. Klee; see also
\cite{N} and \cite{Z}.

\medskip

A great contribution to the study of coverings of
infinite-dimensional Banach spaces was given in 1961 by H. H.
Corson in his classical paper \cite{C}, motivated by topological
reasons. The main result of that paper states that, if a Banach
space $X$ contains some infinite-dimensional reflexive subspace,
then $X$ admits no locally finite covering by bounded convex sets.
Such a result has been recently improved in several directions by
V. P. Fonf and the second author. In fact in \cite{FZ1} it is
proved that, if every compact subset of $X$ meets only finitely
many members of some covering of $X$ by bounded $w$-closed
subsets, then $X$ is $c_0$-saturated. Moreover, in order that
there exist a (algebraically) finite-dimensional compact set that
meets infinitely many members of any covering $\tau$ of $X$ by
closed convex bounded (in short CCB) sets, it is enough for $X$ to
contain an infinite-dimensional separable dual space (see
\cite{FZ2}). In particular, in this case even a segment exists that meets
infinitely many members of $\tau$ whenever the members of $\tau$
are CCB rotund or smooth bodies (see \cite{FZ3}), or simply CCB bodies
if $\tau$ is a tiling and $X$ itself is (infinite-dimensional) reflexive (see \cite {N}).

\medskip

Despite these results, an interesting question asked by Corson
in \cite{C} still remains unanswered. After proving his theorem, he
essentially asks  whether some
infinite-dimensional reflexive space exists that admits a
point-finite covering by CCB bodies. Note that, without any
assumption on the infinite-dimensional Banach space $X$, even
locally finite tilings by CCB bodies can be exhibited, like the
``lattice'' tiling of $c_0$ by suitable translates of the unit
ball (see \cite{F} for a significant characterization of the
separable Banach spaces admitting such tilings as those being
isomorphically polyhedral). Moreover, in special non-separable
spaces even (not locally finite) tilings of order 1 can be found by CCB bodies: see
the surprising construction given by V. Klee in \cite{K2} of a tiling of
$\l_1(\Gamma)$ for suitable $\Gamma$ by pairwise disjoint translates of the
closed unit ball.

\medskip

The aim of this paper is to answer Corson's question in the
affirmative. In fact we prove the following

\bigskip

{\bf Theorem} \ \ {\it Every Banach space $X$ admits a tiling of
order 2 by closed convex bounded bodies.}

\bigskip

Our proof is obtained by combining in a suitable way two main
ideas. The first one allowed V. P. Fonf, A. Pezzotta and the
second author to prove in \cite{FPZ} that any Banach space can be
tiled by CCB bodies. Via the second one, A. H. Stone constructed
in \cite{S} a tiling of order 2 of ${\R}^n$, $n$ any natural
number, by CCB bodies.

\bigskip

Throughout the paper we use standard Geometry of Banach Spaces
notation as in \cite{JLHB}. All the Banach spaces under
consideration are assumed to be real.

\bigskip

\bigskip

\section{Proof of the Theorem}

For finite-dimensional spaces a construction can be found in \cite{S} (the bodies produced there
are also uniformly bounded),
so we can work only in the infinite-dimensional setting.

\bigskip

We recall that, for a normed space $X$ with norm $|| \cdot ||$ and
$0 < \alpha \leq 1$, a set $M \subset S_{X^*}$ is called
$\alpha-$norming if \ ${\rm sup}\{ |f(x)|: f \in M \} \geq \alpha
||x||$ \ for every $x \in X$. A {\it norming set} for $X$ is a
subset of $S_{X^*}$ which is $\alpha-$norming for some $\alpha$.
Moreover, ${\rm norm}(X)$ is the smallest cardinal number $c$ such
that there exists a norming set $M$ for $X$ with $|M| = c$.

Finally, ${\rm dens}(X)$ is the smallest cardinal number $c$ such
that there exists a set $W \subset X$, with $|W| = c$, that is
(strongly) dense in $X$.

It easy to see that \ ${\rm norm}(X) \leq {\rm dens}(X)$ \ for
every Banach space $X$; the equality holds whenever $X$ is weakly
compactly generated (see \cite{L}, Prop. 2.2).

We split our proof into four steps.
\bigskip

\medskip

{\bf Step 1} \ Let us begin producing a special covering
$\widetilde{\sigma}$ by pairwise disjoint convex bounded bodies of $l_\infty
(\Gamma)$, $\Gamma$ any infinite set. Given any family $\{ [a_\nu
,b_\nu ]\} _{\nu \in \Gamma }$ of bounded closed non trivial real
intervals indexed on $\Gamma$ with ${\rm inf}\{b_\nu - a_\nu: \nu
\in \Gamma\} > 0$, we say that the CCB body

\begin{equation}\label{box}
S = \{ t  \in l_\infty (\Gamma):
t(\nu) \in [a_\nu ,b_\nu ],\  \nu \in \Gamma\}
\end{equation}
is a ``box" in $l_\infty (\Gamma).$

Following Stone (see \cite{S}), if in (\ref{box}) exactly for one value of $\nu \in \Gamma$ the
closed interval  $[a_\nu ,b_\nu ]$ is
replaced by the left-open interval  $(a_\nu ,b_\nu ]$, we say that the corresponding set $S$ is a ``lidless box".

We can assume that the set $\Gamma$ is well -ordered, i.e.

$$\Gamma=\{\nu: 1 \leq \nu < \gamma \}, \ \ \ \gamma \ {\rm a \ limit \ ordinal}. $$

Our construction of $\widetilde{\sigma}$ is totally inspired by
the work that was already done in proving Proposition 1.6 in
\cite{FPZ}. Practically, we repeat that construction just replacing
boxes by lidless boxes. It is worthwhile to present it again,
since it will be used in Step 2 and some more work on it will be
done in Step 3. A picture giving an idea of how the construction
works can be found in that paper.
\medskip

Put $A_0 = B_{l_\infty (\Gamma)}$ and, for $\nu \in \Gamma $ and
$n = 0,1,2,...$

\begin{equation}\label{lidtile} A^{(n)}_\nu = \{ t  \in l_\infty (\Gamma): |t(\mu) | \leq 2^n \ {\rm if} \ \mu < \nu,
\ 2^n < t(\nu) \leq 2^{n+1}, \ |t(\mu) | \leq 2^{n+1} \ {\rm if} \ \mu > \nu \}.
\end{equation}

\bigskip

The collection
$$\widetilde{\sigma} = \{A_0, \ \pm A^{(n)}_\nu : \nu \in \Gamma, n =0,1,2,...\}$$
covers $l_\infty (\Gamma)$ and its members are pairwise disjoint.
In fact, trivially any point in $B_{l_\infty (\Gamma)}$ belongs to
no member of $\widetilde{\sigma}$ different from $A_0$. Moreover,
let $ t  \in l_\infty (\Gamma)$ with $||t|| > 1$ and

- in case ${\rm log}_2||t||$ is not an integer, let $\delta$ the
first index in $\Gamma$ such that  $\displaystyle |t(\delta)| >
2^{[{\rm log}_2||t||]}$;

- in case ${\rm log}_2||t||$ is an integer, let $\delta$ the first
index in $\Gamma$ such that $\displaystyle |t(\delta)| > ||t||/2$.

Clearly $t$ belongs to $({\rm sgn}\,  t(\delta))\, A^{([{\rm
log}_2 |t(\delta)|])}_\delta$ $\langle$ resp. to $({\rm sgn}\,
t(\delta))\, A^{({\rm log}_2|t(\delta)|-1)}_\delta \rangle $ if
${\rm log}_2| t(\delta)|$ is not an integer $\langle$ resp.  ${\rm
log}_2|t(\delta)|$ is an integer $\rangle$  and to no other member
of $\widetilde{\sigma}$.

\bigskip

\medskip

{\bf Step 2} \ Let us show how the covering $\widetilde{\sigma}$
that we have built in Step 1 plays a crucial role in providing a
covering $\sigma$ by pairwise disjoint convex bounded bodies for
any Banach space. We need to recall here under the chief headings
what has been done in \cite{FPZ}, Sect. 2, where detailed proofs
are available.

Let $X$ be a normed space. For a suitable $\Gamma$ with $|\Gamma| =  {\rm norm}(X)$, we want
to construct an isomorphic embedding

$$T: X \to  l_\infty (\Gamma)$$
such that the family

$$\{ T^{-1}(A): A \in \widetilde{\sigma}\}$$
provides the desired covering $\sigma$ for $X$.

Let $M$ be a norming set for $X$ with
$|M| = {\rm norm}(X)$; passing to the equivalent norm $|||x||| =
{\rm sup}\{|f(x)|: f \in M\}$, we may assume that $M$ is
$1-$norming. With the aid of Zorn's lemma, for some ordinal $\gamma$ with
$|\gamma| \leq {\rm dens}(X)$ we construct a totally ordered set
of pairs $\{(x_\nu,f_\nu)\}_{1 \leq \nu < \gamma}$ (which, in some
sense, apes a biorthogonal system) such that

(1) $x_\nu \in S(X)$ and $f_\nu \in M, \ 1 \leq \nu < \gamma$;

(2) $|f_\mu(x_\nu)| \leq 1/2, \  1 \leq  \mu < \nu < \gamma$;

(3) $f_\nu(x_\nu) \geq 3/4, \ 1 \leq \nu < \gamma$;

(4) the set $\{f_\nu\}_{1 \leq \nu < \gamma}$ is (1/2)-norming for
$X$.

Set $\Gamma = \{f_\nu: 1 \leq \nu < \gamma\}$. Since $\Gamma$ is a
norming set for $X$, it follows that $|\Gamma| = |M| = {\rm
norm}(X)$. The map \ $T: X \to  l_\infty (\Gamma)$ \ defined as follows

$$(T(x))(f_\nu ) = f_\nu (x), \ \ 1 \leq \nu < \gamma, \ \ x \in X$$
actually is an isomorphic embedding of $X$ into $l_\infty (\Gamma)$,
since we have $(1/2)||x|| \leq ||T(x)||
\leq ||x||$ for every $x \in X$.

Now, let us consider the covering $\widetilde{\sigma}$ of
$l_\infty (\Gamma)$ that has been constructed in Step 1. It is
obvious that $A_0 \cap T(X)$ has nonempty interior relative to
$T(X)$. Moreover, it can be easily seen (see \cite{FPZ}) that, for fixed $1
\leq \nu < \gamma$ and $ n = 0, 1, 2,...$, the point
$$z^{(n)}_\nu = (2^{n+1} -
0.4)T(x_\nu)$$ is an interior point of $A_\nu^{(n)}
\in \widetilde{\sigma}$, so it is an interior point of $A_\nu^{(n)} \cap
T(X)$ relative to $T(X)$ too.

Hence, for every $A \in \widetilde{\sigma}$, the set $T^{-1}(A)$
is a convex bounded body in $X$, so the family
$$\sigma = \{T^{-1}(A): A \in \widetilde{\sigma} \}$$
provides a covering of $X$
by pairwise disjoint convex bounded bodies.

\bigskip

\medskip

{\bf Step 3} \  Following A. H. Stone (see \cite{S}), we now
produce a refinement $\widetilde\tau$ of $\widetilde{\sigma}$ of
order 2, that turns out to be a tiling of $l_\infty (\Gamma)$ by
CCB bodies. To do that, it is enough to express each lidless box $
A^{(n)}_\nu \in \widetilde{\sigma}$ as a countable union of boxes,
in such a way that any point of $ A^{(n)}_\nu$ belongs to at most
two of them.

For $n, \nu$ fixed, let $\epsilon^{(n)}_\nu$ be a positive number
such that

\begin{equation}\label{epsilon}
\epsilon^{(n)}_\nu < 1, \ \ \ z^{(n)}_\nu + 2\epsilon^{(n)}_\nu \,
B_{T(X)} \subset A^{(n)}_\nu.
\end{equation}

Let $\{a^{(n,j)}_\nu\}_{j=0}^\infty$ a strictly decreasing null
sequence of positive numbers such that

$$a^{(n,0)}_\nu = 2^n$$
and

$$z^{(n)}_\nu + \epsilon^{(n)}_\nu B_{T(X)} \subset \{ t  \in A^{(n)}_\nu:
\ 2^n + a^{(n,1)}_\nu < t(\nu) \leq 2^{n+1} \}.$$

For $j=0,1,2,...$ let us set

\begin{equation}\label{tile}
A^{(n,j)}_\nu = \{ t  \in A^{(n)}_\nu: \ 2^n+ a^{(n,j+1)}_\nu \leq
t(\nu) \leq 2^n + a^{(n,j)}_\nu \}.
\end{equation}

Reasoning as in step 1, it is easy to see that the collection of
CCB bodies

$$\widetilde\tau = \{A_0, \ \pm A^{(n,j)}_\nu : \nu \in \Gamma,\ n,j \in {\N} \cup \{0\} \}$$
gives the desired tiling of $l_\infty (\Gamma)$ of order 2.

\medskip

For $\Gamma = \{1,2\}$, the figure gives an idea of how the
construction works.

\begin{figure}[!h]
\begin{center}
\includegraphics[width=12cm]{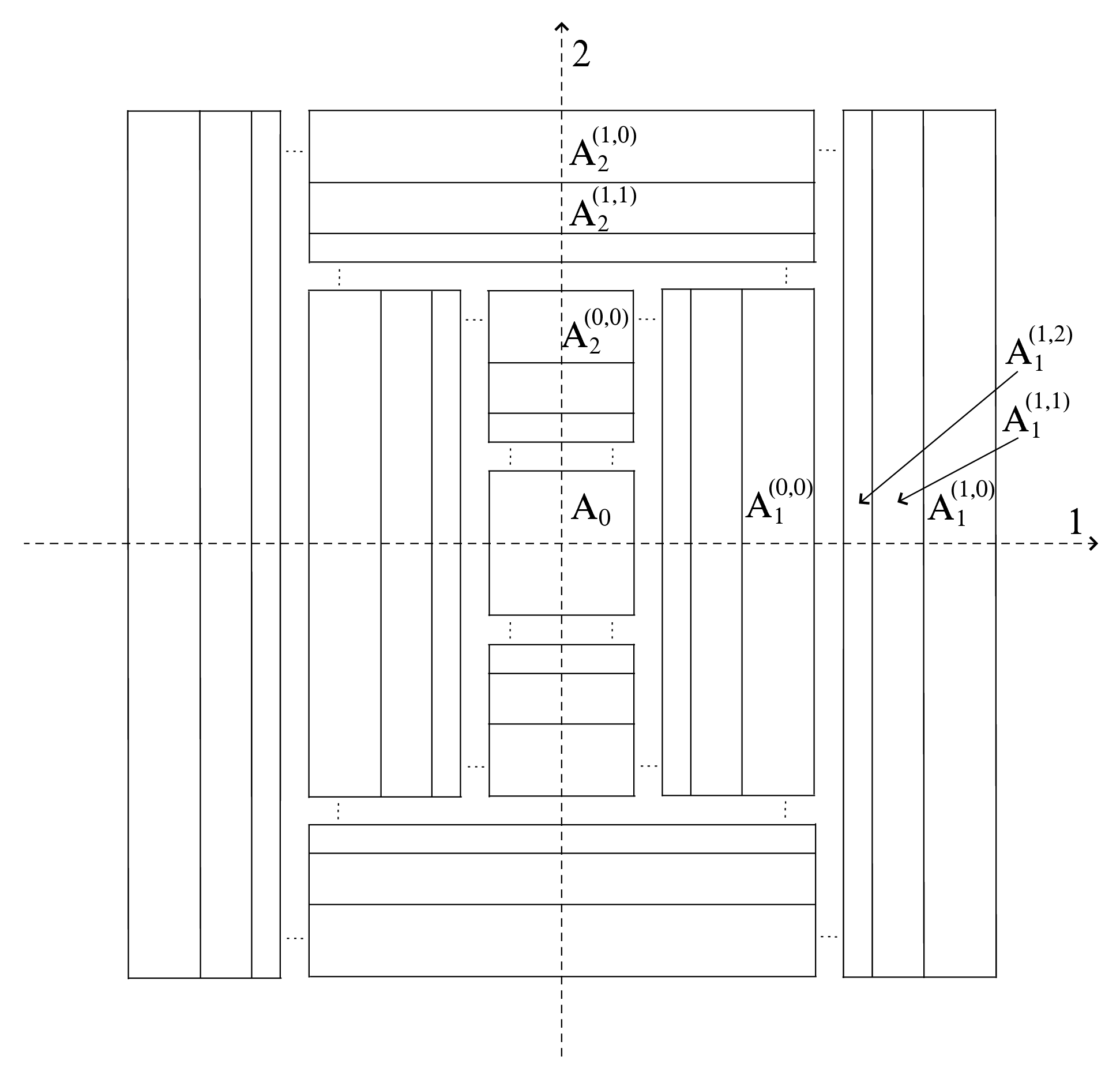}
\end{center}
\end{figure}

\bigskip

\bigskip

{\bf Step 4} \ Finally, we have only to follow the same procedure used
in Step 2, just replacing $\widetilde{\sigma}$ by
$\widetilde{\tau}$ when defining the isomorphic embedding $T$. It
remains only to prove that, for any value of $\nu, n, j$, also the
CCB body $A^{(n,j)}_\nu$ has nonempty interior relative to $T(X)$.
This is clear for $j=0$. For $j=1,2,...$, it is easy to
show that the segment $S$ of $T(X)$, having the origin and
$z^{(n)}_\nu$ as its endpoints, meets the interior of
$A^{(n,j)}_\nu \cap T(X)$ relative to $T(X)$. In fact set

$$r^{(n,j)}_\nu =
{{1} \over {2}}\, (a^{(n,j)}_\nu - a^{(n,j+1)}_\nu)$$ and let
$z^{(n,j)}_\nu$ be the point in $S$ such that

$$z^{(n,j)}_\nu(\nu) = 2^n + a^{(n,j+1)}_\nu + r^{(n,j)}_\nu.$$

From (\ref{lidtile}), (\ref{epsilon}) and (\ref{tile}) it easily
follows that

$$z^{(n,j)}_\nu + {\rm min}\{1, r^{(n,j)}_\nu\}\, \epsilon^{(n)}_\nu B_{T(X)} \subset A^{(n,j)}_\nu.$$

Hence, for every $A \in \widetilde\tau$, the set $T^{-1}(A)$ is a
CCB body in $X$, so the family
$$\tau = \{ T^{-1}(A): A \in
\widetilde\tau \}$$ is a tiling of $X$ of order 2. \ \ {$\Box$}

\bigskip
\bigskip

\bigskip

{\bf Remark} \ \   Clearly our construction cannot lead
anyway to a tiling $\tau$ uniformly bounded from above or
from below, i.e. we cannot get that the members of
$\tau$ are uniformly bounded or that there exists some positive
$r$ such that all of them contain some ball of radius $r$. So,
while the choice of coefficients $2^n$ in \cite{FPZ} was suggested
by the possibility to get all members of $\tau$ uniformly bounded
from below, in our Step 1 it has been made just in order to
refer quickly to that paper for those proofs that have been
omitted here. It is worthwhile to notice that in \cite{FL} (Prop. 2.8) it
is suggested the way to obtain, in any Banach space with the Radon-Nykodym
property, a tiling uniformly bounded from above just applying a cutting procedure
to each member of a given tiling in a straightforward transfinite way. This can obviously
be done on the members of our tiling $\tau$. After that, to each slice that we have
obtained, we can apply a new cutting procedure as described in Step 3, just referred
to the linear continuous functional we used to produce the slice.

\bigskip

\bigskip

Andrea Marchese

Dipartimento di Matematica

Universit\`a degli Studi

Largo B. Pontecorvo, 5

56127 Pisa PI, Italy

E-mail address: marchese@mail.dm.unipi.it

ph: ++39 050 2213 229

\bigskip

\bigskip

Clemente Zanco

Dipartimento di Matematica

Universit\`a degli Studi

Via C. Saldini, 50

20133 Milano MI, Italy

E-mail address: clemente.zanco@unimi.it

ph. ++39 02 503 16164 \ \ \ \ fax  ++39 02 503 16090

\end{document}